\documentclass[11pt, dvipdfmx]{amsart}
\usepackage[utf8]{inputenc}
\usepackage{textcomp}
\usepackage{mathtools, amsmath, amssymb, mathrsfs, bm}
\usepackage{tikz}
\usepackage{hyperref}
\usepackage{amsthm}
\usepackage{kantlipsum} 

\calclayout

\usepackage{color}

\theoremstyle{definition}
\newtheorem{counter}{counter}[section]
\newtheorem{definition}[counter]{Definition}
\newtheorem{theorem}[counter]{Theorem}

\newtheorem{lemma}[counter]{Lemma}

\newtheorem*{Proof}{Proof}

\newtheorem{corollary}[counter]{Corollary}

\makeatletter
\newcommand{\tpitchfork}{%
  \vbox{
    \baselineskip\z@skip
    \lineskip-.52ex
    \lineskiplimit\maxdimen
    \m@th
    \ialign{##\crcr\hidewidth\smash{$-$}\hidewidth\crcr$\pitchfork$\crcr}
  }%
}
\makeatother

\renewcommand{\qed}{\hfill \square}

\renewcommand{\le}{\leqslant}
\renewcommand{\ge}{\geqslant}
\renewcommand{\phi}{\varphi}
\renewcommand{\epsilon}{\varepsilon}

\newcommand{\cal}{\mathcal}

\newcommand{\R}{\mathbb R}
\newcommand{\C}{\mathbb C}

\newcommand{\I}{\textrm{I}}
\newcommand{\II}{\textrm{I\hspace{-1pt}I}}

\title{Unstability problem of real analytic maps}
\author{Karim Bekka}
\address{(K.~Bekka) Institut de recherche math\'ematique de Rennes, Universit\'e de Rennes1, Campus Beaulieu, 35042 Rennes cedex, France}
\email{karim.bekka@univ-rennes1.fr}
\author{Satoshi Koike}
\address{(S.~Koike) Department of Mathematics, Hyogo University of Teacher Education, 942-1 Shimokume, Kato, Hyogo 673-1494, Japan}
\email{koike@hyogo-u.ac.jp}
\author{Toru Ohmoto}
\address{(T.~Ohmoto) Department of Applied Mathematics, Waseda University, 3-4-1 Okubo, Shinjuku-ku, Tokyo 169-8555, Japan}
\email{toruohmoto@waseda.jp}
\author{Masahiro Shiota}
\address{(M.~Shiota) Graduate School of Mathematics, Nagoya University, Furo-cho,
Chigusa-ku, Nagoya 464-8602, Japan}
\author{Masato Tanabe}
\address{(M.~Tanabe) Department of Mathematics, Graduate School of Science, Hokkaido University, Kita 10 Nishi 8, Kita-ku, Sapporo, Hokkaido, 060-0810, Japan}
\email{tanabe.masato.i8@elms.hokudai.ac.jp}
\date{\today}

\begin{document}

\maketitle

\vspace{-12pt}
\begin{center}
\small{{\em To the memory of Masahiro Shiota}}
\end{center}

\begin{abstract}
As well-known, the $C^\infty$ stability of proper $C^\infty$ maps is characterized by the infinitesimal $C^\infty$ stability. In the present paper we study the counterpart in real analytic context. In particular, we show that the infinitesimal $C^\omega$ stability does not imply $C^\omega$ stability; for instance, {\em a Whitney umbrella $\R^2 \to \R^3$ is not $C^\omega$ stable}. 
A main tool for the proof is a relative version of Whitney's Analytic Approximation Theorem which is shown by using H.~Cartan's Theorems A and B. 
\end{abstract}


\section{Introduction}

Singularity theory of smooth maps was initiated by H.~Whitney and R.~Thom in the mid-20th century.
One of the most significant works in this field is structual stability theorem of $C^\infty$ maps established by J.~N.~Mather \cite{MatherI, MatherII, MatherIII, MatherIV, MatherV, MatherVI}.
At the beginning of the theory, the notion of $C^\infty$ stability was most important.  
We say that $C^\infty$ maps $f, g \colon M \to N$ are {\em $\cal{A}$-equivalent} if there exist $C^\infty$ diffeomorphisms $\sigma, \tau$ of $M$ and $N$, respectively, so that $g = \tau \circ f \circ \sigma^{-1}$. 
Then, $f$ is {\em $C^\infty$ stable} if any $C^\infty$ maps $g$ sufficiently close to $f$ in the Whitney $C^\infty$-topology is $\cal{A}$-equivalent to $f$. The $C^\infty$ stability of a proper map is characterized in terms of infinitesimal deformations of the map, namely, {\em the $C^\infty$ stability and the infinitesimal $C^\infty$ stability are equivalent}. Furthermore, the stability is essentially reduced to the infinitesimal stability of multi-germs $f \colon (M, S) \to (N, q)$ for any finite sets $S \subset f^{-1}(q)$ and any $q \in N$.

In the present paper, we study the counterpart in the context of real analytic category; in particular, we are interested in how the rigidity of analytic maps makes difference for the stability. In fact, the $C^\omega$ stability is quite restrictive; we show that {\em the infinitesimal $C^\omega$ stability does not imply $C^\omega$ stability.}


Following works of Kiyoshi Oka and Henri P. Cartan which largely developed sheaf theoretic approach to real and complex analytic geometry, a major principle is that ``locally analytic'' implies ``globally analytic''. 
Besides, real analytic geometry is much more subtle than complex analytic geometry. 
In real case, thanks to Malgrange's preparation theorem (cf.~Mather \cite{MatherIV}), we know as a local property of $C^\omega$ map-germs
\begin{center}
infinitesimal $C^\infty$-stability $\Leftrightarrow$ infinitesimal $C^\omega$-stability 
\end{center}
\cite[Lemma 11.3]{MatherStr}.
However, we will see that the global nature of $C^\omega$ maps is completely different: 
\begin{center}
$C^\infty$-stability $\not\Rightarrow$ $C^\omega$-stability. 
\end{center}
Precisely, the main result of the present paper is stated as follows.
\begin{theorem}\label{main}
\emph{
Let $M$ and $N$ be $C^\omega$ manifolds and $f \colon M \to N$ a proper $C^\omega$ map which is $C^\infty$ stable. 
If $f$ is not an immersion in case of $m < n$, or $f$ has an $A_3$-singular point (swallowtail singularity) in case of $m \ge n$, then $f$ is $C^\omega$ unstable.
}
\end{theorem}
The proof heavily relies on famous \emph{Cartan's Theorems A and B} \cite{Cartan, FS}. 
The most crucial fact is that the image (or critical value set) of a $C^\omega$ map $f$ satisfying the condition in Theorem \ref{main} is {\em semianalytic} but {\em not analytic} in the target manifold. 

As a remark for the complex analytic case, it is reasonable to consider holomorphic maps from Stein manifolds to open domains in $\C^n$ (cf. e.g., \cite{For}), because we need enough many holomorphic functions on manifolds in order to pursue an analogy to the $C^\infty$ stability. In the complex case, the image of the critical locus of a proper holomorphic map is always complex analytic, that has a clear contrast to the nature of real analytic maps, on which we focus in the present paper. 
The semianalyticity is also treated for a different study of $C^\omega$ map-germs in \cite{Damon}. 


\section{Preliminaries}

\subsection{$C^\omega$ Stability}

Let $M$ and $N$ be $C^\omega$ manifolds of dimension $m$ and $n$, respectively. 

\begin{definition}
(1) Two $C^\omega$ maps $f, g \colon M \to N$ are \emph{$C^\omega$-$\cal{A}$-equivalent} if $g = \tau \circ f \circ \sigma^{-1}$ for some $C^\omega$ diffeomorphisms $\sigma$ and $\tau$ of $M$ and $N$, respectively.\\
(2) We say that a $C^\omega$ map $f$ is \emph{$C^\omega$ stable} if any $C^\omega$ maps $g$ sufficiently close to $f$ in the Whitney $C^\infty$-topology is $C^\omega$-$\cal{A}$-equivalent to $f$.
\end{definition}

As mentioned in the introduction, infinitesimal $C^\omega$ stability and infinitesimal $C^\infty$ stability are equivalent for $C^\omega$ map germs.
In particular, an infinitesimal stable $C^\omega$ map-germ admits a polynomial normal form via some $C^\omega$ diffeomorphism-germ of the source and the target, according to the finite determinacy \cite{MatherIII}. Here we recall two typical singularity types -- {\em Whitney umbrella} and {\em Swallowtail}. Indeed, in our proof of Theorem \ref{main},  the main step (\S \ref{proof}) is to show that any representative of those singularity types are not $C^\omega$ stable. 

\vspace{5pt}

\noindent 
-- {\bf  \emph{Whitney umbrella}}: 
This is the simplest stable singularity in case $\kappa:=n-m >0$, denoted by $A_1$ (Thom-Boardman symbol $\Sigma^{1, 0}$); the normal form is given by the $(m-\kappa-1)$-parameter trivial unfolding of 
\[A_1 \colon (x, y_1, \dots, y_{\kappa}) \mapsto (x^2, xy_1, \dots xy_{\kappa}, y_1, \dots, y_{\kappa}).\]
The image is a semianalytic subset, not analytic one. In fact, let $(u, v_1, \cdots, v_\kappa, w_1, \cdots, w_\kappa)$ denote coordinates of the target space of the above normal form, then the analytic closure is given by $v_i^2=uw_i^2\, (1\le i \le \kappa)$, while the actual image of the map lies on the area $u\ge 0$. Namely, the half line of the $u$-axis in the analytic closure is missing in the image of the map -- 
that is the {\em shank} (stick) of the umbrella. 
Note that the analytic closure is not coherent (e.g., J.-C.~Tougeron \cite[Remarques 7.8 of Chapitre II]{Tougeron}, Damon \cite[\S 1]{Damon}).
We also note that a special case of Theorem \ref{main} is shown by M.~Shiota \cite[Fact 3.8]{Shiota}. Its proof is constructive and based on the non-coherency.

\vspace{5pt}

\noindent 
-- {\bf \emph{Swallowtail singularity}}: This is a stable singularity of type $A_3$ (or $\Sigma^{m - n + 1, 1,1, 0}$) in case $m \ge n$; the normal form  is given by the $(n-3)$-parameter trivial unfolding of 
\[A_3 \colon (x, y, z, \bm{w}) \mapsto (x^4 + y x^2 + z x + Q(\bm{w}), y, z),\]
where $\bm{w} = (w_1, \dots, w_{m - n})$ and $Q(\bm{w}) = \pm w_1^2 \pm \dots \pm w_{m - n}^2$. 
The image of the critical locus is the well-known swallowtail surface in $\R^3$;  it is a semianalytic subset, and its analytic closure is given by the discriminant of the quartic equation $x^4+ux^2+vx+w=0$ (taking coordinates $(-w, u, v)$ of $\R^3$). 
This analytic subset contains the double curve defined by $4w=u^2$ and $v=0$, 
and the half of the curve with $u<0$ is out of the swallowtail surface.  
In this case also, the analytic closure is not coherent.

Combining Theorem \ref{main} and the above observations in Whitney umbrella and swallowtail singularity, we have the following corollary on an analogy to Mather's nice range \cite{MatherV}. 

\begin{corollary} 
{\em 
Let $M$ be a compact $C^\omega$ manifold of dimension $m$, and let $N$ be a $C^\omega$ manifold of dimension $n$. If $m$ and $n$ satisfy the condition $3 \le n \le m$ or $2\le m < n < 2m$, then the pair $(m,n)$ of positive integers is not in the nice range of dimensions with respect to the $C^\omega$ stability. 
}
\end{corollary} 

\vspace{5pt}

\subsection{Key facts in the $C^\omega$ category}
For a $C^\omega$ manifold $X$, let $\cal{O}_X$ (resp.~$\cal{E}_X$) denote the sheaf of all germs of $C^\omega$ (resp.~$C^\infty$) functions on $X$.

\begin{theorem}[Whitney's Approximation Theorem \cite{Whitney}] \label{approx-0}
{\em 
For $i = 1, 2$, let $X_i$ be a $C^\omega$ manifold.
Then for any $C^\infty$ map $f \colon X_1 \to X_2$, we can find a $C^\omega$ map $g \colon X_1 \to X_2$ arbitrarily close to $f$ for the Whitney $C^\infty$-topology.
}
\end{theorem}

\begin{theorem}[Cartan's Theorem {\cite{Cartan}, see also \cite[\S 2.1]{FS}}] \label{Cartan}
{\em 
Let $\cal{M}$ be a coherent sheaf of $\cal{O}_X$-modules on a $C^\omega$ manifold $X$.
Then the following hold: 
\begin{itemize}
\item[(A)] at each point $x \in X$, it holds that $\cal{M}_x = H^0(X; \cal{M}) \cdot \cal{O}_{X, x}$;
\item[(B)] $H^1 (X; \cal{M}) = 0$.
\end{itemize}
}
\end{theorem}

On the above theorem, we notice that 
\begin{enumerate}
\item if $X$ is compact, then there exists a generating system of $H^0(X, \cal{M})$ consisting of finite members, for Noetherianness of $H^0(X; \cal{O}_X)$ \cite[Th\'eor\`eme (I,9)]{Frisch};
\item if $\cal{M}$ is an ideal sheaf of $\cal{O}_X$, then the natural map $H^0(X; \cal{O}_X) \to H^0(X; \cal{O}_X / \cal{M})$ is surjective.
\end{enumerate}

As an application of these facts, we also have the following.

\begin{theorem}[Tognoli's Approximation Theorem {\cite[Theorem 3.1]{Tognoli}} (see also {\cite[\S 3]{BKS}})] \label{Tognoli}
{\em 
Let $\cal{M}$ be a coherent sheaf of $\cal{O}_X$-modules on a $C^\omega$ manifold $X$.
For any open set $U \subset X$, it holds that $\Gamma(U; \cal{M})$ is dense in $\Gamma(U; \cal{M} \otimes_{\cal{O}_X} \cal{E}_X)$ for the Whitney $C^\infty$-topology.
}
\end{theorem}

\section{Proof of Theorem \ref{main}}

\subsection{Relative version of Whitney's approximation theorem} 
We prepare the following general lemmata. We show them by using applications of Cartan's Theorems A and B. 
For the proof of Theorem \ref{main}, we will use a relative version of Whitney's approximation theorem which has the form of Lemma \ref{approx-2}.

\begin{lemma}\label{approx-1}
{\em 
For $i = 1, 2$, let $X_i \subset \R^n$ be a compact $C^\omega$ submanifold (with possibly boundary).
Also let $Y_i$ be a regular $C^\omega$ submanifold of $X_i$, which is a closed subset in $X_i$. 
Suppose that there is a $C^\infty$ diffeomorphism $\phi \colon X_1 \to X_2$ whose restriction also induces the $C^\infty$ diffeomorphism $\phi|_{Y_1} \colon Y_1 \to Y_2$. 
Then there is a $C^\omega$ diffeomorphism $\psi \colon X_1 \to X_2$ so that $\psi|_{Y_1} \colon Y_1 \to Y_2$ is $C^\omega$ diffeomorphic.
Moreover, if there is a $C^\omega$ diffeomorphism $\phi^{(1)} \colon Y_1 \to Y_2$ sufficiently close to $\phi|_{Y_1} \colon Y_1 \to Y_2$, then the $C^\omega$ diffeomorphism $\psi \colon X_1 \to X_2$ can be chosen as an extension of $\phi^{(1)}$.
}
\end{lemma}

\begin{Proof}
By assumption, we take a $C^\infty$ diffeomorphism $\phi \colon (X_1, Y_1) \to (X_2, Y_2)$.
We modify $\phi$ so that $\phi|_{Y_1} \colon Y_1 \to Y_2$ is $C^\omega$ diffeomorphic, by going through the following procedure:
\begin{enumerate}
\item Using Theorem \ref{approx-0}, we approximate $\phi|_{Y_1} \colon Y_1 \to Y_2$ by a $C^\omega$ diffeomorphism $\phi^{(1)} \colon Y_1 \to Y_2$; 
\item Using a partition of unity on $X_1$, we extend $\phi^{(1)}$ to a $C^\infty$ diffeomorphism $\phi^{(2)} \colon (X_1, Y_1) \to (X_2, Y_2)$ and then rewrite $\phi \coloneqq \phi^{(2)}$.
\end{enumerate}

Applying Theorem \ref{Cartan} (B) to the ideal sheaf $\cal{I}_{Y_1}$ of $Y_1$ on $X_1$, we extend $\phi|_{Y_1}$ to a $C^\omega$ map $\Phi \colon X_1 \to \R^n$; we obtain a $C^\infty$ map $\Phi - \phi \colon X_1 \to \R^n$ which vanishes on $Y_1$.
Namely, every component function of $\Phi - \phi$ belongs to $\cal{I}_{Y_1} \otimes_{\cal{O}_{X_1}} \cal{E}_{X_1}$.
Using Theorem \ref{Tognoli}, we approximate $\Phi - \phi$ by a $C^\omega$ map $\Psi$ whose components are in $\cal{I}_{Y_1}$.

We now put $\psi \coloneqq \Phi - \Psi \colon X_1 \to \R^n$.
Since $\psi$ is close to $\phi$ on $X_1$, its image $\psi(X_1)$ is contained in a $C^\omega$ tubular neighborhood of $X_2$ in $\R^n$.
Composing $\psi$ to the projection of the neighborhood, we have a $C^\omega$ map between $X_1$ and $X_2$; we also write this map as $\psi$.
Since $\psi$ is close to $\phi$, it is diffeomorphic as well.
Moreover, since $\psi$ coincides with $\phi$ on $Y_1$, it sends $Y_1$ onto $Y_2$.
This completes the proof. $\qed$
\end{Proof}

\begin{lemma}\label{approx-2}
{\em 
For $i = 1, 2$, let $X_i \subset \R^n$ be a $C^\omega$ submanifold with boundary $\partial X_i$. 
Also let $Y_i$ be a regular $C^\omega$ submanifold of $X_i$, which is a closed subset in $X$ and transverse to $\partial X_i$. 
Suppose that there is a $C^\infty$ diffeomorphism $\phi \colon X_1 \to X_2$ whose restrictions also induce $C^\infty$ diffeomorphisms $\phi|_{Y_1} \colon Y_1 \to Y_2$ and $\phi|_{\partial X_1} \colon \partial X_1 \to \partial X_2$. 
Then there is a $C^\omega$ diffeomorphism $\psi \colon X_1 \to X_2$ so that both $\psi|_{Y_1} \colon Y_1 \to Y_2$ and $\psi|_{\partial X_1} \colon \partial X_1 \to \partial X_2$ are $C^\omega$ diffeomorphic.
}
\end{lemma}

\begin{Proof}
For each $i = 1, 2$, put $Z_i = Y_i \cup \partial X_i$ and $W_i = Y_i \cap \partial X_i$. Note that $W_i$ is a regular $C^\omega$ submanifold of $Y_i$, $\partial X_i$, and $X_i$.
By assumption, we take a $C^\infty$ diffeomorphism $\phi \colon (X_1; Y_1, \partial X_1) \to (X_2; Y_2, \partial X_2)$.
We modify $\phi$ so that both $\phi|_{Y_1} \colon Y_1 \to Y_2$ and $\phi|_{\partial X_1} \colon \partial X_1 \to \partial X_2$ are $C^\omega$ diffeomorphic, by going through the following procedure: 
\begin{enumerate}
\item Using Theorem \ref{approx-0}, we approximate $\phi|_{W_1} \colon W_1 \to W_2$ by a $C^\omega$ diffeomorphism $\phi^{(1)} \colon W_1 \to W_2$;
\item We extend $\phi^{(1)}$ to a homeomorphism $\phi^{(2)} \colon (Z_1, W_1) \to (Z_2, W_2)$ which is close to $\phi|_{Z_1}$ and whose restrictions $\phi^{(2)}|_{Y_1} \colon$ $Y_1 \to Y_2$ and $\phi^{(2)}|_{\partial X_1} \colon \partial X_1 \to \partial X_2$ are both $C^\infty$ diffeomorphic;
\item Applying Lemma \ref{approx-1} to both $\phi^{(2)}|_{Y_1} \colon Y_1 \to Y_2$ and $\phi^{(2)}|_{\partial X_1} \colon \partial X_1 \to \partial X_2$, we extend $\phi^{(1)}$ to a homeomorphism $\phi^{(3)} \colon Z_1 \to Z_2$ which is close to $\phi^{(2)}$ and whose restrictions $\phi^{(3)}|_{Y_1} \colon Y_1 \to Y_2$ and $\phi^{(3)}|_{\partial X_1} \colon \partial X_1 \to \partial X_2$ are both $C^\omega$ diffeomorphic;
\item We extend $\phi^{(3)}$ to a $C^\infty$ diffeomorphism $\phi^{(4)} \colon (X_1; Y_1, \partial X_1) \to (X_2; Y_2, \partial X_2)$ and then rewrite $\phi \coloneqq \phi^{(4)}$.
\end{enumerate}

Since $Z_1$ has only normal crossing singularities, we have that the ideal sheaf $\cal{I}_{Z_1}$ of $Z_1$ on $X_1$ is coherent. By the same reason, $\phi|_{Z_1}$ admits an extension to some $C^\omega$ map around each point of $Z_1$.
Then, applying Theorem \ref{Cartan} (B) to the ideal sheaf $\cal{I}_{Z_1}$ of $Z_1$ on $X_1$, we extend $\phi|_{Z_1}$ to a $C^\omega$ map $\Phi \colon X_1 \to \R^n$.
The remained argument is the same as the proof of Lemma \ref{approx-1}. $\qed$
\end{Proof}

\subsection{Necessary condition for $C^\omega$ stability}\label{proof}
We prove the theorem only in the case of $m = 2$ and $n = 3$ because other cases can be proved similarly. 
The essential point is that the critical value set is semianalytic but not analytic. 

Let $f \colon M \to N$ be a proper $C^\omega$ map which is $C^\infty$ stable.
Suppose that $f$ is not an immersion ($M$ is a surface and $N$ is a $3$-fold).
Since $f$ is of class $C^\infty$, it has a singular point, say $p \in M$, at which $f$ is of type Whitney umbrella. 
We show that $f$ is $C^\omega$ unstable. 
More precisely, by performing a real analytic surgery of $f$ on a neighborhood of $p$ in a certain way, 
we show that there exists a $C^\omega$ map $g \colon M \to N$ arbitrarily close to $f$ in
the Whitney $C^\infty$-topology which is not $C^\omega$ equivalent to $f$.

Take analytic coordinate neighborhoods around $p \in U$ and $f(p) \in V$ so that $f|_U \colon U \to V$ is written as $(x, y) \mapsto (u, v, w) = (x^2, xy, y)$. The analytic closure in $V$ of $f(U)$ is given by $uw^2 = v^2$; the u-axis is included, although the half line $u < 0$ of the axis is not in $f(U)$. Take a point $q \in V$ on the half line close to $0$ and an open ball $Q_1 \subset V$ centered at $q$ which does not intersect with the image of $f$. Also take $Q_2 \subset Q_1$ a smaller open ball centered at $q$.
Put $X_1 \coloneqq \overline{Q_1} - Q_2$ and $X_2 \coloneqq \overline{D(1)} - D(1/2)$, where $D(r)$ denotes the standard open $3$-ball centered at $0$ with radius $r$. Both $X_1$ and $X_2$ are compact $C^\omega$ manifolds with boundary.
Let $Y_1 \coloneqq \{(u, 0, 0) \in \overline{Q_1}\}$, a segment, and $Y_2$ an analytic
curve in $D(1)$ having a singular point at $0$ and being transverse to $\partial X_2$.
Then $Y'_i \coloneqq X_i \cap Y_i$ ($i = 1, 2$) consists of two pieces of non-singular $C^\omega$ curves which are transverse to $\partial X_i$. Since $(X_1, Y'_1)$ and $(X_2, Y'_2)$ are $C^\infty$ diffeomorphic, from Lemma \ref{approx-2}, we can find a $C^\omega$ diffeomorphism $\psi : (X_1, Y'_1) \to (X_2, Y'_2)$.

Regarding $\psi$ as a glueing map, we define an abstract $C^\omega$ manifold $N' \coloneqq (N - Q_2) \cup_\psi D(1)$. Note that there is a $C^\infty$ diffeomorphism $\rho \colon N \to N'$ which is identical on $N - Q_1$. Moreover, $N'$ can be embedded in some $\R^\ell$ as a closed $C^\omega$ submanifold, i.e., it is a Stein manifold (Grauert \cite{Grauert}); so we have an injective $C^\infty$ map $\rho \colon N \to \R^\ell$.
By Whitney's approximation theorem, we find a $C^\omega$ approximation $\rho_0 \colon N \to \R^\ell$ sufficiently close to $\rho$.
Although $\rho(N)$ does not coincide with $N'$ in general, $\rho_0(N)$ is contained in a $C^\omega$ tubular neighborhood of $N'$. Composing it with an appropriate orthogonal projection to $N'$, we denote the resulting $C^\omega$ map by $\bar{\rho} \colon N \to N'$.
The map $\bar{\rho}$ is sufficiently close to $\rho$, hence it is $C^\omega$ diffeomorphic. Note that $\bar{\rho}$ is no longer identical on $N - Q_1$ but close to the identity.

Since $f(U)$ is contained in $N - Q$, we have a proper $C^\omega$ map
\[\bar{f} \coloneqq \iota \circ f \colon M \to N - Q_1 \to N',\]
where $\iota \colon N - Q_1 \to N'$ is the inclusion. Obviously, $\bar{f}(U) = f(U) \subset N - Q_1$.
Put $V' \coloneqq (V - Q_2) \cup_\psi D(1) \subset N'$. Remember that the analytic closure of $f(U)$ in $V$ contains the straight line along $Y_1$, while the analytic closure of $\bar{f}(U)$ in $V'$ contains the singular curve $Y_2$ (the analytic closure of $\bar{f}(U)$ contains $X_2 \cap Y_2$ by the identification via $\psi$ with $X_1 \cap Y_1$, and hence it contains the entire curve $Y_2$).

We define $g \coloneqq \bar{\rho}^{-1} \circ \bar{f} \colon M \to N$; it is also a proper $C^\omega$ map which does not coincide with $f$, but sufficiently close to $f$ for the Whitney $C^\infty$-topology.
Now, suppose that $f$ is $C^\omega$ stable. By definition, $g$ is $C^\omega$ equivalent to $f$; there are $C^\omega$ diffeomorphisms $\sigma$ and $\tau$ on $M$ and $N$, respectively, such that $g = \tau \circ f \circ \sigma^{-1}$. Consider their complexification $f_\C, g_\C \colon M_\C \to N_\C$, they are also holomorphically equivalent. However, $(f|_U)_\C$ has only a singular point $p \in U_\C$, while $(g|_{\sigma(U)})_\C$ has another singular point which corresponds via $\bar{\rho}_\C$ to the singular point $0$ of $(Y_2)_\C$ in $V'_\C$, that makes a contradiction. 
This completes the proof. $\qed$

\subsection{Remark on a sufficient condition for $C^\omega$ stability} 
We say that a singularity type $\eta$ of maps is \emph{adjacent} to another type $\tau$ 
if the $\cal{A}$-orbit of type $\eta$ is included in the closure of the $\cal{A}$-orbit of type $\tau$; we denote it by $\eta \to \tau$.
Let $f: M \to N$ be a $C^\omega$ map which is $C^\infty$ stable. 


In case of $m > n$, there are stable singularity types $A_k, D_k, E_k, \cdots$ in $\Sigma^{m - n + 1}$ and types $S_6, S_7, \cdots$ in $\Sigma^{m - n + 2}$, and so on, and they have adjacency relations $A_k \to A_{k + 1}$ and $A_3 \to D_4 \to D_5 \to D_6, E_6, S_6, \cdots$, see \cite{Arnold} for the details. 
Hence, the swallowtail $A_3$ is adjacent to any other stable singularities than $A_1$ and $A_2$. 

In case of $m = n$, stable singularities in $\Sigma^1$ are of type $A_k\, (k\ge 1)$ only. The simplest types in $\Sigma^2$ are the elliptic type ($+$) denoted by $\I_{2,2}$ and the hyperbolic type ($-$) denoted by $\II_{2,2}$ , which have codimension $4$:
\[\I_{2, 2}\;\; \mbox{and}\; \; \II_{2, 2} \colon (x, y, u, v, z) \mapsto (x^2 \pm y^2 + ux + vy, xy, u, v, z).\]
There is a remarkable difference between these two types: 
the swallowtail $A_3$ is adjacent to $\I_{2, 2}$ but not to $\II_{2, 2}$, see e.g., \cite{Arnold}. 
The type  $\I_{2, 2}$ is adjacent to other singularity types of codimension $\ge 5$ in $\Sigma^{\ge 2}$ ($\I_{2, 3}, \I_{2, 4}, \II_{2, 4}, \I_{3, 3}$ etc.), hence $A_3$ is also adjacent to them.
Therefore, only $A_1$, $A_2$ and $\II_{2, 2}$ are singularity types to which the swallowtail $A_3$ is not adjacent.
Moreover, we can check that for each germ of type {\em fold} $A_1$, {\em cusp} $A_2$ or {\em hyperbolic} $\II_{2, 2}$, 
the germ of the discriminant set is analytic by direct computation. 

To determine the $C^\omega$-nice range exactly, we have to check the remaining cases $(m, 1)$, $(m,2)$ and $(m, 2m)$ with  $m\ge 2$. This is a bit technical and will be discussed elsewhere.

\end{document}